\documentclass{amsart}
\usepackage{amscd,amssymb,epsfig}
\usepackage{epic,eepic}

\numberwithin{equation}{subsection}

\newtheorem{propo}{Proposition}[subsection]
\newtheorem{corol}[propo]{Corollary}

\newtheorem{lemma}[propo]{Lemma}
\theoremstyle{definition}
\newtheorem{defin}[propo]{Definition}
\newtheorem{examp}[propo]{Example}
\theoremstyle{remark}
\newtheorem{remar}[propo]{Remark}

\newcounter{figureb}

\newcommand{\trid}{\triangleright}
\newcommand{\G}{\mathcal{G}}
\newcommand{\ZZ}{\mathbb{Z}}
\newcommand{\RR}{\mathbb{R}}
\newcommand{\FF}{\mathbb{F}}

\newcommand{\ttt}{\mathbf{t}}
\newcommand{\Hom}{\operatorname{Hom}}

\newcommand{\lkn}{\operatorname{lk}}
\newcommand{\flecha}[3]{#1 \raisebox{-0.3mm}{$\substack{#2\\ \raisebox{1.2mm}{$\longrightarrow$}}$} #3}
\newcommand{\GG}{\mathfrak{G}}
\newcommand{\invar}[2]{\iint #1\,#2}
\begin{document}
\title{Knot theory for self-indexed graphs}
\author[Gra\~na and Turaev]{Mat\'\i as Gra\~na and Vladimir Turaev}
\thanks{The work of M.G. was supported by CONICET (Argentina)}
\address{%
IRMA--CNRS, 7 rue Ren\'e Descartes \newline
\indent 67084 Strasbourg Cedex \newline
\indent France \vspace*{.4cm}\newline
\indent M.G. permanent address:\newline
\indent Depto de Matem\'atica - FCEyN - Universidad de Buenos Aires \newline
\indent Ciudad Universitaria Pab. I \newline
\indent 1428 Buenos Aires \newline
\indent Argentina}
\begin{abstract}
We introduce and study so-called self-indexed graphs. These are (oriented) finite graphs
endowed with a map from the set of edges to the set of vertices. Such graphs naturally
arise from classical knot and link diagrams.  In fact, the graphs resulting from link
diagrams have an additional structure, an integral flow. We call a self-indexed graph with
integral flow a comte.  The analogy with links allows us to define transformations of
comtes generalizing the Reidemeister moves on link diagrams.  We show that many invariants
of links can be generalized to comtes, most notably the linking number, the Alexander
polynomials, the link group, etc. We also discuss finite type invariants and quandle
cocycle invariants of comtes.
\end{abstract}
\maketitle

\section{Introduction}
\subsection{Outline}
By a self-indexed graph, we mean a finite oriented graph provided with a map from
the set of edges to the set of vertices. In this paper we show that the study of such
graphs is closely related to the study of knots and links in Euclidean 3-space $\RR^3$. This
connection allows us to view self-indexed graphs as a generalization of links and to apply
to them a number of well understood tools of knot theory.

Knot theory studies smooth embeddings of $n=1,2,...$ copies of the circle $S^1$
into $\RR^3$.  Two such embeddings are \emph{isotopic\/} if they can be smoothly
deformed into each other in the class of embeddings. An isotopy class of such embeddings is
called an \emph{oriented $n$-component link\/} in $\RR^3$.  (The orientation of the link is
determined by the counterclockwise orientation on $S^1$.) One-component links are called
\emph{knots\/}.  Knots and links are usually presented by their generic projections
to the plane $\RR^2$ where one should keep track of under/overcrossings and the circle
orientations, cf. Figure 5. Such projections are called \emph{oriented link diagrams\/}.
We refer to the monographs \cite{bz,rf} for more on knots and links.

We show in this paper that every oriented link diagram gives rise to a self-indexed graph
with an integral flow. Recall that an \emph{integral flow\/} on an oriented graph is an
assignment of an integer to each arrow such that the algebraic sum of the integers incident
to any vertex is zero.  We call a self-indexed graph with a flow a \emph{comte\/}. The
word comte is the French word for count (nobleman). We use it in our context in analogy
with Russian and German, where the word \emph{graph} has two meanings: the standard one used in
the graph theory and a count (nobleman).  Since \emph{count} has a precise meaning in
mathematics, we use its French translation.

We introduce certain transformations of comtes generalizing the Reidemeister moves on
link diagrams. Since every two diagrams of an oriented link are related by the Reidemeister
moves, we conclude that every oriented link in $\RR^3$ determines a comte, at least up
to our moves. As we shall see, this mapping from the set of links to the set of comtes
(modulo the moves) is essentially injective.  However, the class of comtes is much wider
than the class of links. The theory of comtes can be viewed as a larger, combinatorial
paradigm for knot theory.

We shall show that many standard notions of knot theory extend to comtes. This inlcudes
the link group, the link quandle, the state sum quandle cocycle invariants, the linking
number, the Alexander polynomials, etc.  The case of quandle cocycle invariants of comtes is
treated in more detail. We generalize quandles to so-called \emph{self-indexed q-graphs\/},
define a homology theory for them, and use it to present state sums on comtes as a sort
of double integral. These state sums are invariant under some of our moves, and they are
invariant under all of the moves if the q-graph comes from a quandle.  Indeed, this paper
began by the observation that quandle cocycle invariants could be defined for self-indexed
graphs with flows.

We also briefly discuss finite type invariants of comtes and show that a virtual
link in the sense of Kauffman \cite{kvn} gives rise to a comte.  (However, the resulting
map from the set of virtual links to the set of comtes is not injective.)

\subsection{Definitions}\label{sn:defs}

By a \emph{graph} we mean a $4$-tuple $(V,E,s,t)$, where $V$ and $E$ are finite sets
(resp. of vertices and arrows (=edges)), and $s,t:E\to V$ are maps (the source and the
target). Note that we allow loops (i.e., edges $e$ with $s(e)=t(e)$) and multiple edges
(i.e., edges with the same endpoints).  Each graph gives rise to a finite 1-dimensional
cellular complex, called its \emph{topological realization}.  Its $0$-cells are the vertices
of the graph and its 1-cells are numerated by the arrows of the graph so that the 1-cell
corresponding to an arrow $e$ is an oriented interval leading from $s(e)$ to $t(e)$. As
usual, we shall pictorially present graphs by their topological realizations.

A graph $(V,E,s,t)$ endowed with a map $\ell:E\to V$ is said to be \emph{self-indexed}.
Note that we impose no conditions on the map $\ell$, for instance one can take $\ell=s$
or $\ell=t$, or let $\ell$ be a constant map.  The value of $\ell$ on an edge is called
the \emph{label} of this edge.  We draw edges of a self-indexed graph as on the left-hand
side of the following figure, where $a$ stands for the label of the edge.
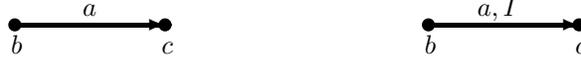
\begin{figure}[ht]\setcounter{figure}{-2}
	\setlength{\unitlength}{0.5cm}
	\begin{picture}(20,1.5)(0,.25)
	\allinethickness{.5mm}
	\put(1,1){\vector(1,0){4}}
	\put(1,1){\circle*{.25}}
	\put(5,1){\circle*{.25}}
	\put(0.9,.25){$b$}
	\put(4.9,.25){$c$}
	\put(2.8,1.25){$a$}
	\put(12,1){\vector(1,0){4}}
	\put(12,1){\circle*{.25}}
	\put(16,1){\circle*{.25}}
	\put(11.9,.25){$b$}
	\put(15.9,.25){$c$}
	\put(13.3,1.25){$a,I$}
	\end{picture}
\end{figure}

An (integral) \emph{flow\/} on a graph $(V,E,s,t)$ is a map $I:E\to\ZZ$ such that for
each vertex $v\in V$, the sum of outgoing flows is equal to the sum of incoming flows:
\begin{equation}\label{eq:cur}
\sum_{e:s(e)=v}I(e)=\sum_{e:t(e)=v}I(e).
\end{equation}
A self-indexed graph with a flow is called a \emph{comte}.  We draw flows as on the right-hand
side of the previous figure, where $I$ stands for the value of the flow on the edge.

Given two self-indexed graphs $\G=(V,E,s,t,\ell)$ and $\G'=(V',E',s',t',\ell')$,
a \emph{homomorphism\/}
$f:\G\to\G'$ is a pair of maps $f_V:V\to V'$, $f_E:E\to E'$
commuting with $s,t,\ell$. More precisely, the following diagrams should commute:
$$\begin{CD}
E @>f_E>> E' &\qquad& E @>f_E>> E' &\qquad& E @>f_E>> E' \\
@VsVV @VVs'V @VtVV @VVt'V @V\ell VV @VV\ell'V \\
V @>f_V>> V', &\qquad& V @>f_V>> V', &\qquad& V @>f_V>> V'.
\end{CD}$$

\bigskip
Clearly, self-indexed graphs and their homomorphisms form a category.

\subsection{Moves on self-indexed graphs and comtes}

We define several transformations of comtes also called \emph{moves}.  In all these
transformations vertices may coincide, but arrows referred to (and drawn) as different
must be different.
\begin{itemize}
\item[R0] A vertex of valence $1$ which is not the label of any arrow can be deleted
	together with the incident arrow, see Figure~\ref{fg:r0}.
	Observe that necessarily the flow of the arrow is $0$.
\item[R1] An arrow labeled by its source or its target can be contracted, see
	Figure~\ref{fg:r1}. The source and target are identified. An arrow pointing from a
	vertex to itself and labeled by this vertex can be deleted.
\item[R2] Two arrows with the same label and the same source (resp. target)
	can be replaced by one arrow as in Figure~\ref{fg:r2a} (resp. Figure~\ref{fg:r2b}).
	The targets (resp. sources) of the arrows are identified and the flows are added.
\item[R3] In presence of an arrow with label $a$, source $b$ and target $t$, any of the
	four arrows in a square with sides labeled $b,a,t,a$ can be removed, see Figure~\ref{fg:r3}.
	The flow is modified as in the figure.
\end{itemize}

We say that two comtes are \emph{isotopic} if they can be related by a sequence of
isomorphisms, moves R0--R3 and the inverse moves. By abuse of language, we shall use the
same word \emph{comte} for an isotopy class of comtes.

We can consider the same transformations for self-indexed graphs, just by ignoring the
flows. We say that two self-indexed graphs are \emph{isotopic} if they can be related by
a sequence of isomorphisms, moves R0--R3 (with flows forgotten) and the inverse moves.

\begin{figure}[ht]\setcounter{figure}{-1}
	\setlength{\unitlength}{0.4cm}
	\begin{picture}(28,10)(-1,0)
	\texture{0101010110101010}
	\allinethickness{0mm}
	\shade\path(0.5,6)(0,6.5)(0,9.5)(0.5,10)(7.5,10)(8,9.5)(8,6.5)(7.5,6)(0.5,6)
	\allinethickness{.5mm}
	\put(0.5,6){\line(1,0){7}}
	\put(0.5,10){\line(1,0){7}}
	\put(0,6.5){\line(0,1){3}}
	\put(8,6.5){\line(0,1){.5}}
	\put(8,8){\line(0,1){1.5}}
	\put(0.5,6.5){\shade\arc{1}{1.57}{3.14}}
	\put(0.5,9.5){\shade\arc{1}{3.14}{4.71}}
	\put(7.5,9.5){\shade\arc{1}{4.71}{6.28}}
	\put(7.5,6.5){\shade\arc{1}{0.00}{1.57}}
	\put(1,8){\makebox(0,0)[l]{\shortstack{A comte $\G$\\without arrows\\labeled by $b$}}}
	\put(8,8){\circle*{.25}}
	\put(8,8){\vector(1,0){3}}
	\put(11,8){\circle*{.25}}
	\put(9,8.25){$t,0$}
	\put(7.85,7.25){$a$}
	\put(10.85,7.25){$b$}
	\put(14,7.75){$\rightleftharpoons$}
	\allinethickness{0mm}
	\shade\path(17.5,6)(17,6.5)(17,9.5)(17.5,10)(24.5,10)(25,9.5)(25,6.5)(24.5,6)(17.5,6)
	\allinethickness{.5mm}
	\put(17.5,6){\line(1,0){7}}
	\put(17.5,10){\line(1,0){7}}
	\put(17,6.5){\line(0,1){3}}
	\put(25,6.5){\line(0,1){.5}}
	\put(25,8){\line(0,1){1.5}}
	\put(17.5,6.5){\shade\arc{1}{1.57}{3.14}}
	\put(17.5,9.5){\shade\arc{1}{3.14}{4.71}}
	\put(24.5,9.5){\shade\arc{1}{4.71}{6.28}}
	\put(24.5,6.5){\shade\arc{1}{0.00}{1.57}}
	\put(24.85,7.25){$a$}
	\put(25,8){\circle*{.25}}
	\put(20.75,8){\makebox(0,0)[l]{\shortstack{$\G$}}}
	\allinethickness{0mm}
	\shade\path(0.0,0)(-0.5,0.5)(-0.5,3.5)(0.0,4)(7.0,4)(7.5,3.5)(7.5,0.5)(7.0,0)(0.0,0)
	\allinethickness{.5mm}
	\put(0.0,0){\line(1,0){7}}
	\put(0.0,4){\line(1,0){7}}
	\put(-0.5,.5){\line(0,1){3}}
	\put(7.5,.5){\line(0,1){.5}}
	\put(7.5,2){\line(0,1){1.5}}
	\put(0.0,0.5){\shade\arc{1}{1.57}{3.14}}
	\put(0.0,3.5){\shade\arc{1}{3.14}{4.71}}
	\put(7.0,3.5){\shade\arc{1}{4.71}{6.28}}
	\put(7.0,0.5){\shade\arc{1}{0.00}{1.57}}
	\put(0.5,2){\makebox(-0.5,0)[l]{\shortstack{A comte $\G$\\without arrows\\labeled by $b$}}}
	\put(7.5,2){\circle*{.25}}
	\put(10.5,2){\vector(-1,0){3}}
	\put(10.5,2){\circle*{.25}}
	\put(8.75,2.25){$t,0$}
	\put(7.35,1.25){$a$}
	\put(10.35,1.25){$b$}
	\put(13.5,1.75){$\rightleftharpoons$}
	\allinethickness{0mm}
	\shade\path(17.0,0)(16.5,.5)(16.5,3.5)(17.0,4)(24.0,4)(24.5,3.5)(24.5,0.5)(24.0,0)(17.0,0)
	\allinethickness{.5mm}
	\put(17.0,0){\line(1,0){7}}
	\put(17.0,4){\line(1,0){7}}
	\put(16.5,0.5){\line(0,1){3}}
	\put(24.5,0.5){\line(0,1){.5}}
	\put(24.5,2){\line(0,1){1.5}}
	\put(17.0,0.5){\shade\arc{1}{1.57}{3.14}}
	\put(17.0,3.5){\shade\arc{1}{3.14}{4.71}}
	\put(24.0,3.5){\shade\arc{1}{4.71}{6.28}}
	\put(24.0,0.5){\shade\arc{1}{0.00}{1.57}}
	\put(24.35,1.25){$a$}
	\put(24.5,2){\circle*{.25}}
	\put(20.25,2){\makebox(-0.5,0)[l]{\shortstack{$\G$}}}
\end{picture}
\caption{Move R0}\label{fg:r0}
\end{figure}
\renewcommand{\thefigure}{\arabic{figure}}

\begin{figure}[ht]
	\setlength{\unitlength}{0.5cm}
	\begin{picture}(18,5.5)(-1,3.5)
	\allinethickness{.5mm}
	\put(1,7){\vector(1,0){3}}
	\put(1,7){\circle*{.25}}
	\put(0.8,6){$a$}
	\put(4,7){\circle*{.25}}
	\put(3.8,6){$b$}
	\put(2,7.25){$a,I$}
	\put(6.5,6.75){$\rightleftharpoons$}
	\put(9,7){\circle*{.25}}
	\put(8.2,6){$a=b$}
	\put(11,6.75){$\rightleftharpoons$}
	\put(13.5,7){\vector(1,0){3}}
	\put(13.5,7){\circle*{.25}}
	\put(13.3,6){$a$}
	\put(16.5,7){\circle*{.25}}
	\put(16.3,6){$b$}
	\put(14.5,7.25){$b,I$}
	\put(1.5,4){\circle*{.25}}
	\put(0.8,3.85){$a$}
	\put(4,3.85){$a,I$}
	\spline(1.5,4)(3,5.2)(4,4)(3,2.8)(1.5,4)
	\put(1.65,3.85){\vector(-3,2){.1}}
	\put(6.5,3.85){$\rightleftharpoons$}
	\put(8.8,3){$a$}
	\put(9,4){\circle*{.25}}
\end{picture}
\caption{Move R1}\label{fg:r1}
\end{figure}

\renewcommand{\thefigure}{\arabic{figure}(\textup{\alph{figureb}})}
\begin{figure}[ht]\setcounter{figureb}{1}
	\setlength{\unitlength}{0.5cm}
	\begin{picture}(16,4)(-1,1)
	\allinethickness{.5mm}
	\put(-.25,1.85){$a$}
	\put(0.5,2){\circle*{.25}}
	\put(0.5,2){\vector(2,1){2.5}}
	\put(0.5,2){\vector(2,-1){2.5}}
	\put(3,3.25){\circle*{.25}}
	\put(3.5,3.1){$s$}
	\put(3,0.75){\circle*{.25}}
	\put(3.5,0.6){$t$}
	\put(0.8,3.0){$b,I_1$}
	\put(0.8,0.7){$b,I_2$}
	\put(6,1.85){$\rightleftharpoons$}
	\put(9,2){\circle*{.25}}
	\put(8.85,1){$a$}
	\put(9,2){\vector(1,0){4}}
	\put(13,2){\circle*{.25}}
	\put(12.25,1){$s=t$}
	\put(9.5,2.5){$b,I_1+I_2$}
\end{picture}
\caption{Move R2(a)}\label{fg:r2a}
\end{figure}

\begin{figure}[ht]\addtocounter{figure}{-1}\addtocounter{figureb}{1}
	\setlength{\unitlength}{0.5cm}
	\begin{picture}(16,4)(-1,1)
	\allinethickness{.5mm}
	\put(-.25,1.85){$a$}
	\put(0.5,2){\circle*{.25}}
	\put(3,3.25){\vector(-2,-1){2.5}}
	\put(3,0.75){\vector(-2,1){2.5}}
	\put(3,3.25){\circle*{.25}}
	\put(3.5,3.1){$t$}
	\put(3,0.75){\circle*{.25}}
	\put(3.5,0.6){$s$}
	\put(0.8,3.0){$b,I_1$}
	\put(0.8,0.7){$b,I_2$}
	\put(6,1.85){$\rightleftharpoons$}
	\put(9,2){\circle*{.25}}
	\put(8.85,1){$a$}
	\put(13,2){\vector(-1,0){4}}
	\put(13,2){\circle*{.25}}
	\put(12.25,1){$s=t$}
	\put(9.7,2.5){$b,I_1+I_2$}
\end{picture}
\caption{Move R2(b)}\label{fg:r2b}
\end{figure}
\renewcommand{\thefigure}{\arabic{figure}}

\setcounter{figure}{2}
\begin{figure}[ht]
\setlength{\unitlength}{0.5cm}
\begin{picture}(24.5,27)(-3,0)
	\allinethickness{.5mm}
	\put(6.5,26){\circle*{.25}}
	\put(6.5,22){\circle*{.25}}
	\put(6.5,26){\vector(0,-1){4}}
	\put(6.35,26.25){$b$}
	\put(6.35,21.25){$t$}
	\put(6.75,23.85){$a,I_0$}
	\put(-1,22){\circle*{.25}}
	\put(3,22){\circle*{.25}}
	\put(-1,26){\circle*{.25}}
	\put(3,26){\circle*{.25}}
	\put(-1.5,26.25){$c$}
	\put(-1.5,21.5){$u$}
	\put(3.25,21.5){$s$}
	\put(3.25,26.25){$r$}
	\put(-1,22){\vector(1,0){4}}
	\put(0.5,21.25){$t,I_2$}
	\put(-1,26){\vector(0,-1){4}}
	\put(-2.35,23.85){$a,I_1$}
	\put(-1,26){\vector(1,0){4}}
	\put(0.5,26.25){$b,I_3$}
	\put(3,26){\vector(0,-1){4}}
	\put(3.25,23.85){$a,I_4$}
	\put(9,23.35){$\rightleftharpoons$}
	\put(20.5,26){\circle*{.25}}
	\put(20.5,22){\circle*{.25}}
	\put(20.5,26){\vector(0,-1){4}}
	\put(20.35,26.25){$b$}
	\put(20.35,21.25){$t$}
	\put(20.75,23.85){$a,I_0$}
	\put(13,22){\circle*{.25}}
	\put(17,22){\circle*{.25}}
	\put(13,26){\circle*{.25}}
	\put(17,26){\circle*{.25}}
	\put(12.5,26.25){$c$}
	\put(12.5,21.5){$u$}
	\put(17.25,21.5){$s$}
	\put(17.25,26.25){$r$}
	\put(13,22){\vector(1,0){4}}
	\put(13.5,21.25){$t,I_2+I_4$}
	\put(13,26){\vector(0,-1){4}}
	\put(10.15,23.85){$a,I_1+I_4$}
	\put(13,26){\vector(1,0){4}}
	\put(13.5,26.25){$b,I_3-I_4$}
	\allinethickness{.2mm}
	\put(-3,20.5){\line(1,0){25.5}}
	\allinethickness{.5mm}
	\put(9,16.35){$\rightleftharpoons$}
	\put(20.5,19){\circle*{.25}}
	\put(20.5,15){\circle*{.25}}
	\put(20.5,19){\vector(0,-1){4}}
	\put(20.35,19.25){$b$}
	\put(20.35,14.25){$t$}
	\put(20.75,16.85){$a,I_0$}
	\put(13,15){\circle*{.25}}
	\put(17,15){\circle*{.25}}
	\put(13,19){\circle*{.25}}
	\put(17,19){\circle*{.25}}
	\put(12.5,19.25){$c$}
	\put(12.5,14.5){$u$}
	\put(17.25,14.5){$s$}
	\put(17.25,19.25){$r$}
	\put(13,15){\vector(1,0){4}}
	\put(13.5,14.25){$t,I_2+I_3$}
	\put(13,19){\vector(0,-1){4}}
	\put(10.15,16.85){$a,I_1+I_3$}
	\put(17,19){\vector(0,-1){4}}
	\put(17.25,16.75){$a,I_4-I_3$}
	\allinethickness{.2mm}
	\put(-3,13.5){\line(1,0){25.5}}
	\allinethickness{.5mm}
	\put(9,9.35){$\rightleftharpoons$}
	\put(20.5,12){\circle*{.25}}
	\put(20.5,8){\circle*{.25}}
	\put(20.5,12){\vector(0,-1){4}}
	\put(20.35,12.25){$b$}
	\put(20.35,7.25){$t$}
	\put(20.75,9.85){$a,I_0$}
	\put(13,8){\circle*{.25}}
	\put(17,8){\circle*{.25}}
	\put(13,12){\circle*{.25}}
	\put(17,12){\circle*{.25}}
	\put(12.5,12.25){$c$}
	\put(12.5,7.5){$u$}
	\put(17.25,7.5){$s$}
	\put(17.25,12.25){$r$}
	\put(13,8){\vector(1,0){4}}
	\put(13.5,7.25){$t,I_2-I_1$}
	\put(17,12){\vector(0,-1){4}}
	\put(17.25,9.85){$a,I_4+I_1$}
	\put(13,12){\vector(1,0){4}}
	\put(13.5,12.25){$b,I_3+I_1$}
	\allinethickness{.2mm}
	\put(-3,6.5){\line(1,0){25.5}}
	\allinethickness{.5mm}
	\put(9,2.35){$\rightleftharpoons$}
	\put(20.5,5){\circle*{.25}}
	\put(20.5,1){\circle*{.25}}
	\put(20.5,5){\vector(0,-1){4}}
	\put(20.35,5.25){$b$}
	\put(20.35,0.25){$t$}
	\put(20.75,2.85){$a,I_0$}
	\put(13,1){\circle*{.25}}
	\put(17,1){\circle*{.25}}
	\put(13,5){\circle*{.25}}
	\put(17,5){\circle*{.25}}
	\put(12.5,5.25){$c$}
	\put(12.5,0.5){$u$}
	\put(17.25,0.5){$s$}
	\put(17.25,5.25){$r$}
	\put(17,5){\vector(0,-1){4}}
	\put(17.25,2.75){$a,I_4+I_2$}
	\put(13,5){\vector(0,-1){4}}
	\put(10.15,2.85){$a,I_1-I_2$}
	\put(13,5){\vector(1,0){4}}
	\put(13.5,5.25){$b,I_3+I_2$}
\end{picture}
\caption{Move R3}\label{fg:r3}
\end{figure}

\subsection{Comte of a link}
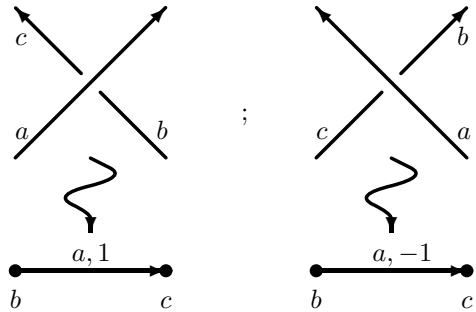
\begin{figure}[ht]
\setlength{\unitlength}{0.5cm}
\begin{picture}(13,8.5)(0,0.5)
	\allinethickness{.5mm}
	\put(1,1){\circle*{.25}}
	\put(5,1){\circle*{.25}}
	\put(1,1){\vector(1,0){4}}
	\put(0.85,0){$b$}
	\put(4.85,0){$c$}
	\put(2.5,1.25){$a,1$}
	\put(1,4){\vector(1,1){4}}
	\put(5,4){\line(-1,1){1.75}}
	\put(2.75,6.25){\vector(-1,1){1.75}}
	\put(1,4.5){$a$}
	\put(4.75,4.5){$b$}
	\put(1,7){$c$}
	\spline(3,4)(4,3.5)(2,3)(3,2.5)(3,2.25)
	\put(3,2.25){\vector(0,-1){.25}}
	\put(7,5){$;$}
	\put(9,1){\circle*{.25}}
	\put(13,1){\circle*{.25}}
	\put(9,1){\vector(1,0){4}}
	\put(8.85,0){$b$}
	\put(12.85,0){$c$}
	\put(10.5,1.25){$a,-1$}
	\put(13,4){\vector(-1,1){4}}
	\put(9,4){\line(1,1){1.75}}
	\put(11.25,6.25){\vector(1,1){1.75}}
	\put(9,4.5){$c$}
	\put(12.75,4.5){$a$}
	\put(12.75,7){$b$}
	\spline(11,4)(12,3.5)(10,3)(11,2.5)(11,2.25)
	\put(11,2.25){\vector(0,-1){.25}}
\end{picture}
\caption{How to pass from a link to a comte}\label{fg:dr5}
\end{figure}

We now explain how every oriented link diagram in $\RR^2$ gives rise to a comte.
Viewed as a subset of the plane, the diagram consists of a finite number of disjoint
embedded oriented arcs.  The set $V$ of these arcs will be the set of vertices of our comte.
Each crossing of the diagram gives rise to three arcs $a,b,c\in V$: the arc $a$ contains
the overpass, the arc $b$ contains the underpass lying on the right of $a$, and $c$
contains the underpass lying on the left of $a$, see Figure \ref{fg:dr5}.  Consider an
arrow from $b$ to $c$ labeled by $a$ with a flow $+1$ if the crossing is positive and $-1$
if the crossing is negative. (Note that the direction of this arrow from $b$ to $c$ is
induced by the orientation of the arcs $b,c$ if the crossing is positive, and is reversed
for a negative crossing.) The set $V$ with such arrows corresponding to all crossings of
the diagram is the comte determined by the diagram.

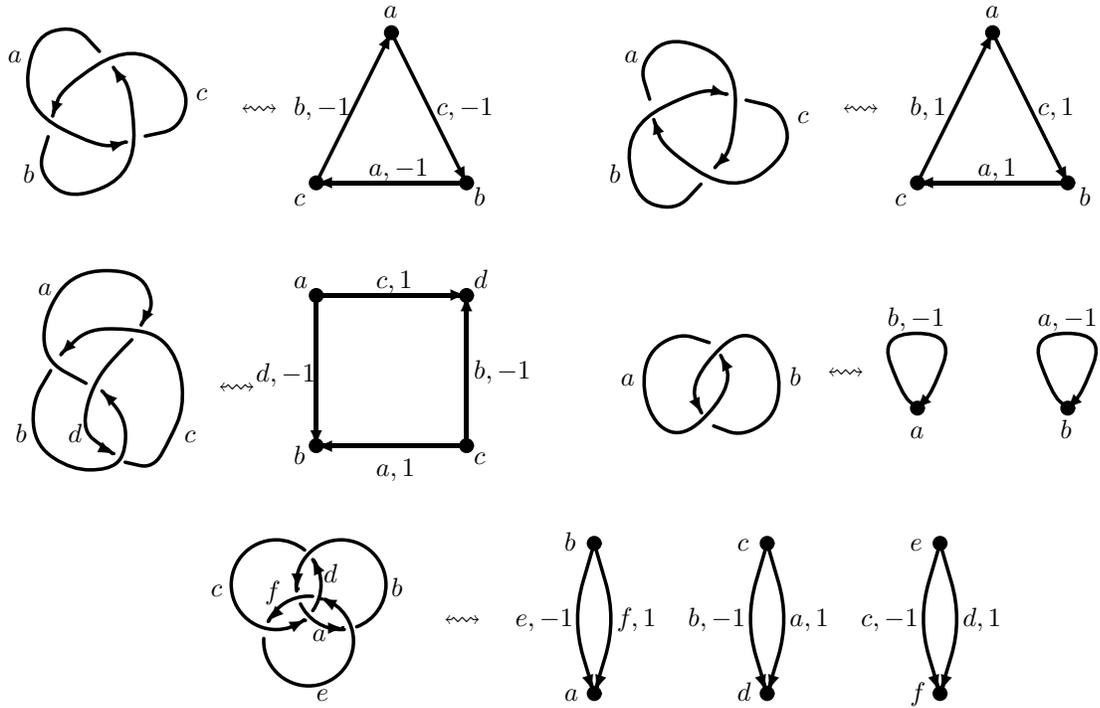
\begin{figure}[ht]
\setlength{\unitlength}{0.01cm}
\allinethickness{.5mm}
\begin{picture}(1500,950)(0,370)
\spline(44,1141)(29,1094)(74,1054)(161,1092)(156,1206)(132,1231)
	\put(132,1231){\vector(-2,3){2}}
	\put(10,1080){$b$}
\spline(112,1255)(79,1291)(22,1274)(12,1178)(112,1126)(146,1134)
	\put(146,1134){\vector(4,1){2}}
	\put(-10,1240){$a$}
\spline(173,1142)(221,1153)(233,1211)(156,1268)(61,1207)(51,1174)
	\put(51,1174){\vector(-1,-3){2}}
	\put(240,1190){$c$}
\put(300,1170){$\leftrightsquigarrow$}
\put(400,1080){\circle*{15}}\put(370,1050){$c$}
	\put(400,1080){\vector(1,2){100}}\put(370,1170){$b,-1$}
\put(600,1080){\circle*{15}}\put(610,1050){$b$}
	\put(600,1080){\vector(-1,0){200}}\put(470,1090){$a,-1$}
\put(500,1280){\circle*{15}}\put(490,1300){$a$}
	\put(500,1280){\vector(1,-2){100}}\put(560,1170){$c,-1$}
\spline(844,1191)(829,1238)(873,1278)(961,1240)(955,1127)(932,1102)
	\put(932,1102){\vector(-1,-1){2}}
	\put(790,1080){$b$}
\spline(912,1078)(878,1041)(822,1059)(812,1155)(912,1206)(945,1198)
	\put(945,1198){\vector(4,-1){2}}
	\put(810,1240){$a$}
\spline(972,1190)(1020,1179)(1033,1121)(956,1064)(861,1126)(851,1159)
	\put(851,1159){\vector(-1,3){2}}
	\put(1040,1160){$c$}
\put(1100,1170){$\leftrightsquigarrow$}
\put(1200,1080){\circle*{15}}\put(1170,1050){$c$}
	\put(1200,1080){\vector(1,2){100}}\put(1190,1170){$b,1$}
\put(1400,1080){\circle*{15}}\put(1415,1050){$b$}
	\put(1400,1080){\vector(-1,0){200}}\put(1280,1090){$a,1$}
\put(1300,1280){\circle*{15}}\put(1290,1300){$a$}
	\put(1300,1280){\vector(1,-2){100}}\put(1360,1170){$c,1$}
\spline(94,813)(38,844)(38,906)(79,964)(164,964)(186,920)(169,894)
	\put(169,894){\vector(-1,-2){2}}
	\put(30,930){$a$}
\spline(47,830)(21,786)(27,728)(76,696)(123,696)(146,713)(149,763)(117,801)
	\put(117,801){\vector(-1,1){2}}
	\put(0,735){$b$}
\spline(146,708)(184,699)(207,740)(225,781)(221,836)(195,877)(149,888)(94,886)(62,853)
	\put(62,853){\vector(-1,-1){2}}
	\put(225,735){$c$}
\spline(155,871)(137,853)(111,824)(94,786)(90,748)(114,728)(129,720)
	\put(129,720){\vector(2,-1){2}}
	\put(70,735){$d$}
\put(270,800){$\leftrightsquigarrow$}
\put(400,930){\circle*{15}}\put(370,940){$a$}
	\put(400,930){\vector(1,0){200}}\put(480,940){$c,1$}
\put(600,930){\circle*{15}}\put(610,940){$d$}
	\put(600,730){\vector(0,1){200}}\put(610,820){$b,-1$}
\put(600,730){\circle*{15}}\put(610,705){$c$}
	\put(600,730){\vector(-1,0){200}}\put(480,690){$a,1$}
\put(400,730){\circle*{15}}\put(370,705){$b$}
	\put(400,930){\vector(0,-1){200}}\put(320,815){$d,-1$}
\spline(911,777)(899,803)(917,841)(960,882)(998,870)(1019,821)(1010,768)(966,742)(929,756)
	\put(911,777){\vector(1,-2){2}}
	\put(1030,810){$b$}
\spline(940,847)(952,821)(934,783)(890,742)(853,754)(833,803)(841,856)(885,882)(923,867)
	\put(940,847){\vector(-1,2){2}}
	\put(805,810){$a$}
\put(1080,820){$\leftrightsquigarrow$}
\put(1200,780){\circle*{15}}\put(1190,740){$a$}
	\spline(1200,780)(1220,800)(1250,880)(1150,880)(1180,800)(1200,780)
	\put(1200,780){\vector(-1,-1){2}}\put(1160,890){$b,-1$}
\put(1400,780){\circle*{15}}\put(1390,740){$b$}
	\spline(1400,780)(1420,800)(1450,880)(1350,880)(1380,800)(1400,780)
	\put(1400,780){\vector(-1,-1){2}}\put(1360,890){$a,-1$}
\put(433,545){\arc{120}{1.5}{2.7}}
	\put(435,487){\vector(1,0){2}}
	\put(395,470){$a$}
\put(433,545){\arc{120}{3.0}{7.5}}
	\put(375,540){\vector(0,-1){2}}
	\put(500,530){$b$}
\put(347,545){\arc{120}{0.9}{5.4}}
	\put(384,497){\vector(2,1){2}}
	\put(260,530){$c$}
\put(347,545){\arc{120}{5.7}{6.9}}
	\put(397,575){\vector(-1,2){2}}
	\put(410,548){$d$}
\put(390,470){\arc{120}{-1.2}{3.2}}
	\put(338,498){\vector(-1,-1){2}}
	\put(400,390){$e$}
\put(390,470){\arc{120}{3.6}{4.8}}
	\put(408,526){\vector(-2,1){2}}
	\put(332,525){$f$}
\put(570,490){$\leftrightsquigarrow$}
\put(770,400){\circle*{15}}\put(730,390){$a$}
\put(770,600){\circle*{15}}\put(730,590){$b$}
	\spline(770,400)(740,500)(770,600)
	\spline(770,400)(800,500)(770,600)
	\put(766,410){\vector(1,-3){2}}\put(665,490){$e,-1$}
	\put(774,410){\vector(-1,-3){2}}\put(800,490){$f,1$}
\put(1000,400){\circle*{15}}\put(960,390){$d$}
\put(1000,600){\circle*{15}}\put(960,590){$c$}
	\spline(1000,400)(970,500)(1000,600)
	\spline(1000,400)(1030,500)(1000,600)
	\put(996,410){\vector(1,-3){2}}\put(895,490){$b,-1$}
	\put(1004,410){\vector(-1,-3){2}}\put(1030,490){$a,1$}
\put(1230,400){\circle*{15}}\put(1190,390){$f$}
\put(1230,600){\circle*{15}}\put(1190,590){$e$}
	\spline(1230,400)(1200,500)(1230,600)
	\spline(1230,400)(1260,500)(1230,600)
	\put(1226,410){\vector(1,-3){2}}\put(1125,490){$c,-1$}
	\put(1234,410){\vector(-1,-3){2}}\put(1260,490){$d,1$}
\end{picture}
\caption{Several link diagrams and their comtes}\label{fg:fl}
\end{figure}

We draw in Figure~\ref{fg:fl} the comtes associated to several simple knot and link diagrams.
As an exercise the reader may verify that the topological realization of the comte determined
by a diagram of an $n$-component link is a disjoint union of $n$ circles.

It is easy to observe that two oriented link diagrams in $\RR^2$ presenting isotopic
oriented links in $\RR^3$ give rise to isotopic comtes.  Indeed, such diagrams are related
by a finite sequence of oriented Reidemeister moves.  We need only to prove that under
these moves the comte changes by a sequence of our moves R0--R3 and the inverse moves.
The action of the first (resp. second, third) Reidemeister move on the comte can be
achieved with R1 (resp.  R2 + R0, R3 + R0).  We check it for the second Reidemeister move
with both strands oriented in the same direction in Fig.~\ref{fg:srr2} and leave the other
cases to the reader.  Thus the isotopy class of the comte derived from an oriented link
diagram depends only on the (isotopy class of the) link itself.  In this way we obtain a
map from the set of (isotopy classes of) oriented links in $\RR^3$ to the set of isotopy
classes of comtes. Forgetting the flows, we obtain a map from the set of oriented links
in $\RR^3$ to the set of isotopy classes of self-indexed graphs.

\begin{figure}[ht]
\setlength{\unitlength}{0.1cm}
\allinethickness{.5mm}
\begin{picture}(130,24)(0,0)
	\spline(10,0)(20,10)(10,20)\put(10,20){\vector(-1,1){2}}
	\path(20,0)(17,3)\spline(13,7)(10,10)(13,13)\path(17,17)(20,20)\put(20,20){\vector(1,1){2}}
	\put(7,0){$a$}\put(21,0){$b$}\put(7,9){$c$}\put(21,18){$d$}
	\put(30,9){$\rightleftharpoons$}
	\path(40,0)(40,20)
	\path(50,0)(50,20)
	\put(40,20){\vector(0,1){2}}\put(50,20){\vector(0,1){2}}
	\put(37,0){$a$}\put(51,0){$b\!=\!d$}
	\put(60,9){$\leftrightsquigarrow$}
	\put(73,10){\circle*{1}}\put(70,9){$c$}
	\put(79,16){\vector(-1,-1){6}}\put(79,16){\circle*{1}}\put(80,14){$d$}\put(69,14){$a,-1$}
	\put(79,4){\vector(-1,1){6}}\put(79,4){\circle*{1}}\put(80,2){$b$}\put(72,3){$a,1$}
	\put(83,9){$\rightleftharpoons$}
	\put(93,10){\circle*{1}}\put(92,7){$c$}
	\put(103,10){\circle*{1}}\put(100,7){$b\!=\!d$}
	\put(103,10){\vector(-1,0){10}}\put(96,11){$a,0$}
	\put(113,9){$\rightleftharpoons$}
	\put(123,10){\circle*{1}}\put(120,7){$b\!=\!d$}
\end{picture}
\caption{Second Reidemeister move $\leftrightsquigarrow$ R2+R0}\label{fg:srr2}
\end{figure}

\begin{remar}
The move R3 can be split as a composition of moves R3(a) and R3(b) below. It is sometimes
easier to work with R3(a) and R3(b) rather than with R3.
\begin{itemize}
\item[R3(a)] In presence of an arrow with label $a$, source $b$ and target $t$,
	we remove any of the four arrows, with a $0$ flow, in a square with sides
	labeled $b,a,t,a$, see Figure~\ref{fg:r3a} (where
	the relation is depicted for one of the sides; analogous figures should be drawn
	for the other three sides).
\item[R3(b)] In presence of an arrow with label $a$, source $b$ and target $t$, we
	shift the flow in a square with sides labeled $b,a,t,a $, see Figure~\ref{fg:r3b}.
\end{itemize}
\end{remar}

\renewcommand{\thefigure}{\arabic{figure}(\textup{3\alph{figureb}})}
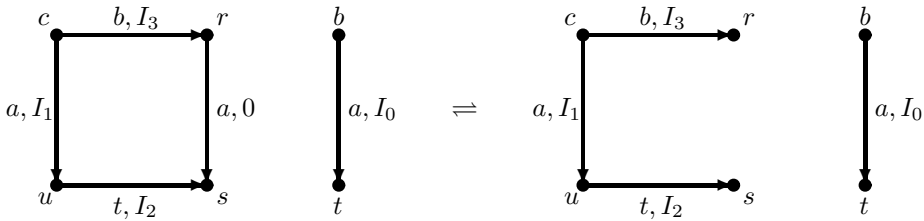
\begin{figure}[ht]\setcounter{figureb}{1}
\setlength{\unitlength}{0.5cm}
\begin{picture}(26,6)(-1,3.5)
	\allinethickness{.5mm}
	\put(8.5,8){\circle*{.25}}
	\put(8.5,4){\circle*{.25}}
	\put(8.5,8){\vector(0,-1){4}}
	\put(8.35,8.25){$b$}
	\put(8.35,3.25){$t$}
	\put(8.75,5.85){$a,I_0$}
	\put(1,4){\circle*{.25}}
	\put(5,4){\circle*{.25}}
	\put(1,8){\circle*{.25}}
	\put(5,8){\circle*{.25}}
	\put(0.5,8.25){$c$}
	\put(0.5,3.5){$u$}
	\put(5.25,3.5){$s$}
	\put(5.25,8.25){$r$}
	\put(1,4){\vector(1,0){4}}
	\put(2.5,3.25){$t,I_2$}
	\put(1,8){\vector(0,-1){4}}
	\put(-0.35,5.85){$a,I_1$}
	\put(1,8){\vector(1,0){4}}
	\put(2.5,8.25){$b,I_3$}
	\put(5,8){\vector(0,-1){4}}
	\put(5.25,5.85){$a,0$}
	\put(11.5,5.85){$\rightleftharpoons$}
	\put(22.5,8){\circle*{.25}}
	\put(22.5,4){\circle*{.25}}
	\put(22.5,8){\vector(0,-1){4}}
	\put(22.35,8.25){$b$}
	\put(22.35,3.25){$t$}
	\put(22.75,5.85){$a,I_0$}
	\put(15,4){\circle*{.25}}
	\put(19,4){\circle*{.25}}
	\put(15,8){\circle*{.25}}
	\put(19,8){\circle*{.25}}
	\put(14.5,8.25){$c$}
	\put(14.5,3.5){$u$}
	\put(19.25,3.5){$s$}
	\put(19.25,8.25){$r$}
	\put(15,4){\vector(1,0){4}}
	\put(16.5,3.25){$t,I_2$}
	\put(15,8){\vector(0,-1){4}}
	\put(13.65,5.85){$a,I_1$}
	\put(15,8){\vector(1,0){4}}
	\put(16.5,8.25){$b,I_3$}
\end{picture}
\caption{Relation R3(a) for one of the sides}\label{fg:r3a}
\end{figure}

\begin{figure}[ht]\addtocounter{figure}{-1}\addtocounter{figureb}{1}
\setlength{\unitlength}{0.5cm}
\begin{picture}(28,6)(-1,3.5)
	\allinethickness{.5mm}
	\put(8.5,8){\circle*{.25}}
	\put(8.5,4){\circle*{.25}}
	\put(8.5,8){\vector(0,-1){4}}
	\put(8.35,8.25){$b$}
	\put(8.35,3.25){$t$}
	\put(8.75,5.85){$a,I_0$}
	\put(1,4){\circle*{.25}}
	\put(5,4){\circle*{.25}}
	\put(1,8){\circle*{.25}}
	\put(5,8){\circle*{.25}}
	\put(0.5,8.25){$c$}
	\put(0.5,3.5){$u$}
	\put(5.25,3.5){$s$}
	\put(5.25,8.25){$r$}
	\put(1,4){\vector(1,0){4}}
	\put(2.5,3.25){$t,I_2$}
	\put(1,8){\vector(0,-1){4}}
	\put(-.35,5.85){$a,I_1$}
	\put(1,8){\vector(1,0){4}}
	\put(2.5,8.25){$b,I_3$}
	\put(5,8){\vector(0,-1){4}}
	\put(5.25,5.85){$a,I_4$}
	\put(12,5.85){$\rightleftharpoons$}
	\put(24.5,8){\circle*{.25}}
	\put(24.5,4){\circle*{.25}}
	\put(24.5,8){\vector(0,-1){4}}
	\put(24.35,8.25){$b$}
	\put(24.35,3.25){$t$}
	\put(24.75,5.85){$a,I_0$}
	\put(17,4){\circle*{.25}}
	\put(21,4){\circle*{.25}}
	\put(17,8){\circle*{.25}}
	\put(21,8){\circle*{.25}}
	\put(16.5,8.25){$c$}
	\put(16.5,3.5){$u$}
	\put(21.25,3.5){$s$}
	\put(21.25,8.25){$r$}
	\put(17,4){\vector(1,0){4}}
	\put(17.5,3.25){$t,I_2+J$}
	\put(17,8){\vector(0,-1){4}}
	\put(14.15,5.85){$a,I_1+J$}
	\put(17,8){\vector(1,0){4}}
	\put(17.5,8.25){$b,I_3-J$}
	\put(21,8){\vector(0,-1){4}}
	\put(21.25,5.85){$a,I_4-J$}
\end{picture}
\caption{Relation R3(b)}\label{fg:r3b}
\end{figure}
\renewcommand{\thefigure}{\arabic{figure}}

\section{Isotopy invariants of self-indexed graphs and comtes}
In this section we generalize a number of well-known invariants of links to self-indexed
graphs and comtes.  We begin with invariants which do not depend on flows.

\subsection{The group of a self-indexed graph}\label{sn:gaq}
Let $\G=(V,E,s,t,\ell)$ be a self-indexed graph. We define the \emph{group} of $\G$ to be
the group generated by elements of $V$ modulo the relations $ab=ca$ whenever there is an arrow
$\flecha bac$ in $\G$.
\begin{lemma}\label{lm:gii}
The group of $\G$ is preserved under the moves R0--R3 (with flows forgotten).
\end{lemma}
\begin{proof}
For R0, simply notice that the group for the LHS has one more generator, $b$, and one
more relation, $b=tat^{-1}$. Hence one can drop $b$ from the set of generators, as on the
RHS. The arrows in R1 give relations of the type $aa=ba$, i.e., $a=b$, and therefore $b$
can be dropped from the set of generators. The loop in the second line in Fig.~\ref{fg:r1}
contributes the relation $aa=aa$ which is tautological.  In move R2(a) we have on both
sides the same relations $s=bab^{-1}=t$. Move R2(b) is analogous. In move R3, on the first
line of RHS of Fig.~\ref{fg:r3a} we have $u=aca^{-1}$,
$s=tut^{-1}$, $r=bcb^{-1}$, $t=aba^{-1}$. Hence
$$s=taca^{-1}t^{-1}=aba^{-1}aca^{-1}ab^{-1}a^{-1} \\
	=abcb^{-1}a^{-1}=ara^{-1},$$
which is precisely the relation on the LHS determined by the forth arrow of the square.
For the other sides of the square the computations are analogous.
\end{proof}

Abelianizing the group of a self-indexed graph $\G=(V,E,s,t,\ell)$ we obtain a free abelian
group with (free) generators bijectively corresponding to the components of $\G$.  Here by
\emph{components\/} of $\G$ we mean the equivalence classes of the equivalence relation in $V$
generated by $b\sim c$ whenever there is an arrow $\flecha b{\phantom{a}}c$ in $\G$. It is clear that
the components of $\G$ bijectively correspond to the components of its topological realization
$\vert\G\vert$.  Thus, the abelianized group of $\G$ is nothing but $H_0(\vert\G\vert;\ZZ)$.
Clearly, the number of components of $\G$ is invariant under the moves R0--R3.

Note finally that the group of a link in $\RR^3$ is isomorphic to the group of its
self-indexed graph.

\subsection{The quandle of a self-indexed graph}
We show that each self-indexed graph gives rise to a quandle.  We first recall the relevant
definitions.

A \emph{rack\/} is a pair $(X,\trid)$ where $X$ is a set and $\trid:X\times X\to X$ is a
binary operation such that
\begin{itemize}
\item the function $x\trid ?:X\to X$ is bijective for all $x\in X$, and
\item $x\trid(y\trid z)=(x\trid y)\trid(x\trid z)$ for all $x,y,z\in X$.
\end{itemize}
A rack $(X,\trid)$ is a \emph{quandle\/} if
\begin{itemize}
\item $x\trid x=x$ for all $x\in X$.
\end{itemize}

For a set $V$, the \emph{free quandle} of $V$ coincides with the union of the conjugacy
classes of elements of $V$ inside the free group generated by $V$. The functor which assigns
the free quandle to a set is left adjoint to the forgetful functor from quandles to sets.

\medskip
The quandle of a self-indexed graph $\G=(V,E,s,t,\ell)$ is the free quandle of $V$
quotiented out by the relations $a\trid b=c$ for each arrow $\flecha bac$ in $\G$.
\begin{lemma}
The quandle of $\G$ is preserved under the moves R0--R3.
\end{lemma}
\begin{proof}
The same as that for the group of $\G$.
\end{proof}
The quandle of the self-indexed graph of an oriented link coincides with the quandle of
the link, as defined in \cite{j} and \cite{matv}.

As an application of quandles, we show that the self-indexed graph is an almost complete
invariant of an oriented link.  For an oriented link $L$, we denote by $\overline L$
its mirror image with reversed orientation on all components. It is easy to see that
the self-indexed graphs of $L$ and $\overline L$ coincide.  (It suffices to present $L$
by a diagram and to consider its mirror image with respect to a plane orthogonal to the
plane of the diagram.) This shows that the self-indexed graph can not distinguish $L$ from
$\overline L$.  However, this is the only source of links with the same self-indexed graphs.
To state the relevant result, we call two oriented links $L_1,L_2$ \emph{weakly isotopic\/}
if $L_1$ is isotopic to $L_2$ or to $\overline L_2$. A link is \emph{splittable} if it is
a disjoint union of two non-empty links. For instance, all knots are non-splittable.

\begin{corol}\label{co:igi}
If two non-splittable oriented links $L_1,L_2$ in $\RR^3$ have isotopic self-indexed graphs,
then $L_1,L_2$ are weakly isotopic.
\end{corol}
\begin{proof}

Suppose first that $L_1,L_2$ are knots.  Since the quandles of $L_1,L_2$ are determined by
their self-indexed graphs, these quandles are isomorphic to each other.  But the quandle
is a full invariant of oriented knots up to weak isotopy (cf. \cite{j}, \cite {matv}).
Therefore $L_1$ is weakly isotopic to $L_2$.  The same argument works for non-splittable
links.
\end{proof}

We don't know whether the comte always distinguishes $L$ from $\overline L$.  Certain
invariants of comtes suggest that this may be the case, cf. \ref{eq:qci} below.  Note also
that a link is splittable if and only if its self-indexed graph is a disjoint union of
two non-empty self-indexed graphs.

\subsection{The Alexander module and Alexander polynomials of a self-indexed graph}
Any finitely generated group $\pi$ gives rise to a ${\ZZ}[H_1(\pi)]$-module called its
Alexander module. It can be computed by the Fox calculus from any presentation of $\pi$
by generators and relations. This module gives rise to a sequence of elements of the group
ring ${\ZZ}[H_1(\pi)/\text{Tors}]$ called the Alexander polynomials of $\pi$.  In particular, the
group of a self-indexed graph yields a module and a sequence of Alexander polynomials. We
give here a direct definition of these module and polynomials.

Let $\G=(V,E,s,t,\ell)$ be a self-indexed graph.  Denote by $\Lambda$ the ring of integer
Laurent polynomials $\ZZ[\ttt,\ttt^{-1}]$.  The \emph{Alexander module\/} $A(\G)$ of $\G$
is the $\Lambda$-module generated by $V$ modulo the relations $c=\ttt b+(1-\ttt)a$ for
each arrow $\flecha bac\in\G$.

\begin{lemma}
The module $A(\G)$ is invariant under the moves R0--R3 (with flows forgotten) on $\G$.
\end{lemma}
\begin{proof}
Let us check, for instance, that the Alexander module is invariant under R3 (we do it for
one of the sides, the other three being analogous).  On the RHS of Figure~\ref{fg:r3a}
we have the relations
$$
r=\ttt c+(1-\ttt)b,\quad u=\ttt c+(1-\ttt)a,\quad s=\ttt u+(1-\ttt)t,\quad t=\ttt b+(1-\ttt)a.
$$
But then
\begin{align*}
s	&= \ttt^2c+\ttt(1-\ttt)a+\ttt(1-\ttt)b+(1-\ttt)^2a
		=\ttt^2c+\ttt(1-\ttt)b+(1-\ttt)a \\
	&= \ttt r+(1-\ttt)a,
\end{align*}
which is exactly the fifth relation on the LHS.
\end{proof}
We now define the $i$-th Alexander polynomial $\Delta_i(\G)$ of $\G$ for any $i=0,1,2,..$.
Present $A(\G)$ by $\# (V)$ generators and $\#(E)$ relations as above (where $\#(\text{a set})$
is the number of elements of the set). Consider the corresponding
$(\# (E)\times \# (V))$-matrix over $\Lambda$. Let $\Delta_i(\G)\in \Lambda$ be the
greatest common divisor of all minors of rank $\#(V)-i$ of this matrix.
By convention, if $\# (V)-i\leq 0$, then $\Delta_i(\G)=1$;
if $\# (V)-i \geq \# (E)+1$, then $\Delta_i(\G)=0$. Clearly, $\Delta_i(\G)$ is defined up
to multiplication by monomials $\pm t^k$ with $k\in \ZZ$ and $\Delta_{i+1}(\G)$ divides
$\Delta_i(\G)$ for all $i$.  These polynomials are preserved under the moves R0--R3 on
$\G$. For the self-indexed graph of an oriented link $L$, this sequence of polynomials
coincides with the (1-variable) Alexander polynomials of $L$.

\begin{examp}
Consider the self-indexed graph
$$\begin{CD}
\G=(c @<b<< a @>c>> b).
\end{CD}
$$
The relations $c=\ttt a+(1-\ttt)b$ and $b=\ttt a+(1-\ttt)c$ give the matrix of relations
$$\left(\begin{array}{ccc} \ttt & -1 & 1-\ttt \\ \ttt & 1-\ttt & -1 \end{array}\right).$$
Its minors of size $2$ are $\ttt(1-\ttt)+\ttt=-\ttt(\ttt-2)$.  Thus, $\Delta_1(\G)=\ttt-2$. As
a consequence, we see that $\G$ is not isotopic to the self-indexed graph of a knot. Indeed,
the Alexander polynomials of knots (and links) are invariant under the conjugation
$\ttt\mapsto\ttt^{-1}$.
\end{examp}

Substituting $\ttt=1$ one easily obtains that for any self-indexed graph $\G$, the sum of
coefficients of $\Delta_1(\G)$ is 0 if $\G$ has $\geq 2$ components and is $\pm 1$ if $\G$
is connected.

We can similarly define the Alexander module (and Alexander polynomials) of $\G$ with
$n$ variables where $n$ is the number of components of $\G$. Namely, let us enumerate
these components by $1,2,\ldots,n$. The (multi-variable) Alexander module is the
$\ZZ[\ttt_1^{\pm},\ldots,\ttt_n^{\pm}]$-module generated by $V$ modulo the relations
$c=\ttt_ib-(1-\ttt_j)a$ for each arrow $\flecha bac$ where $a$ belongs to the $i$-th
component and $b,c$ belong to the $j$-th component.  It is again straightforward to check
that this module is invariant under the moves R0--R3.

\subsection{Linking numbers}
We now introduce our first invariant of comtes depending on the flow.  We say that an
arrow of a self-indexed graph belongs to a certain component of this graph if its source
(and then its target) belongs to this component.  Note that arrows with the same label
may belong to different components.

Let $\G=(V,E,s,t,\ell)$ be a comte with flow $I$ and components $L_1,\ldots,L_n$ (i.e.,
$V=\cup_{i=1}^nL_i$).  For $1\le i,j\le n$, $i\neq j$, the \emph{linking number $\lkn_{ij}$
of $L_i$ with $L_j$,\/} is the sum of flows of arrows belonging to $L_j$ with label in
$L_i$. In other words,
$$
\lkn_{ij}=\sum_{\substack{a\in E\\s(a)\in L_j\\ \ell(a)\in L_i}}I(a).
$$
It is straightforward to check that the linking numbers are invariant under the moves R0--R3.
For the comte of a link we recover the usual linking numbers. We warn, however, that
linking numbers for comtes need not be symmetric, i.e., in general $\lkn_{ij}\neq \lkn_{ji}$.

\begin{examp}
The comte in Fig. \ref{fg:lkn} has three components $L_1=\{a,b\}$, $L_2=\{c\}$, and
$L_3=\{d\}$. The linking numbers are zero, except for $\lkn_{21}=2$.
\begin{figure}[ht]
	\setlength{\unitlength}{0.4cm}
	\begin{picture}(10,2.5)(-1,-0.5)
	\allinethickness{.5mm}
	\put(0,0){\circle*{.25}}
	\put(0,0.2){\vector(1,0){4}}
	\put(4,-0.2){\vector(-1,0){4}}
	\put(4,0){\circle*{.25}}
	\put(-0.8,0){$a$}
	\put(4.4,0){$b$}
	\put(1.5,0.6){$c,1$}
	\put(1.5,-1.2){$c,1$}
	\put(6,0){\circle*{.25}}
	\put(6.4,0){$c$}
	\put(8,0){\circle*{.25}}
	\put(8.4,0){$d$}
	\end{picture}
\caption{A comte with non-symmetric linking numbers}\label{fg:lkn}
\end{figure}
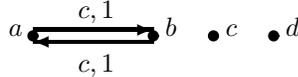
\end{examp}

\subsection{Quandle cocycle invariants of comtes}
For a quandle $X$, a quandle $2$-cocycle of $X$ with values in an abelian group $A$
(written multiplicatively) is a function $f:X\times X\to A$ such that
$$f(x\trid y,x\trid z)f(x,z)=f(x,y\trid z)f(y,z)$$
and $f(x,x)=1$ for all $x,y,z\in X$.  We show that each such $f$ gives rise to an isotopy
invariant of comtes.

Let $\G=(V,E,s,t,\ell)$ be a self-indexed graph with flow $I$.  A \emph{coloring} of $\G$
by $X$ is a function $C:V\to X$ such that for each arrow $\flecha bac$ in $\G$ we have
$C(a)\trid C(b)=C(c)$.  Assume that $X$ is finite. Set
\begin{equation}\label{eq:qci}
\Phi(\G,I,X,f)=\sum_{\substack{\text{colorings}\\C}}\;\prod_{\substack{\text{arrows}\\ \flecha b{a,I}c}}
	f(C(a),C(b))^I \in \ZZ A,
\end{equation}
where the product is taken in $A$ (or in the ring $\ZZ A$), while the addition is that
in $\ZZ A$.  This is a refinement of the invariant counting the colorings of $\G$ by
$X$. Indeed, if we take the map $\epsilon:\ZZ A\to\ZZ$, $\epsilon(\sum_{g\in A} n_gg)=\sum n_g$,
then $\epsilon(\Phi(\G,I,X,f))$ is the number of such colorings.

It is straightforward to see that $\Phi(\G,I,X,f)$ is invariant under moves R0--R3.
By the condition \eqref{eq:cur}, this invariant does not depend on the choice of $f$
in its cohomology class in $H^2(X;A)$.  We shall generalize this invariant in \S\ref{ss:vgd}.

For the comte $(\G,I)$ of an oriented link $L$, the invariant $\Phi(\G,I,X,f)$ coincides
with the invariant $\Phi_f(L)$ defined in \cite{cjkls}.

\begin{remar}
Any rack (in particular a quandle) $X$ yields a self-indexed graph $(X,X\times X,s,t,\ell)$,
where $s(x,y)=y$, $t(x,y)=x\trid y$ and $\ell(x,y)=x$ for any $x,y\in X$.  It is clear that
a coloring of a self-indexed graph $\G$ by a quandle $X$ is nothing but a homomorphism
from $\G$ to the self-indexed graph determined by $X$ (it is also the same as a quandle
homomorphism from the quandle of $\G$ to $X$).
\end{remar}

\section{Virtual links and finite type invariants}
\subsection{The comte of a virtual link}
\begin{figure}[ht]
\setlength{\unitlength}{0.5cm}
\begin{picture}(18,6)(-4,-2.5)
	\allinethickness{.5mm}
	\put(0,0){\circle{6}}
	\put(0,3){\vector(-1,0){.25}}
	\put(-3,0){\vector(2,1){4.8}}
	\put(0,-3){\vector(-1,2){2.4}}
	\put(3,0){\vector(-2,-1){4.8}}
	\put(-0.2,1.7){$+$}
	\put(-1.2,-0.3){$-$}
	\put(1.9,-0.2){$+$}
	\put(-0.3,3.3){$a$}
	\put(-3.5,-0.3){$b$}
	\put(2.4,-2.4){$c$}
	\put(6,0){$\leftrightsquigarrow$}
	\put(9,-2){\circle*{.25}}
	\put(11,2){\circle*{.25}}
	\put(13,-2){\circle*{.25}}
	\put(8.4,-2.6){$a$}
	\put(10.85,2.25){$b$}
	\put(13.2,-2.6){$c$}
	\put(11,2){\vector(-1,-2){2}}
	\put(11,2){\vector(1,-2){2}}
	\put(13,-2){\vector(-1,0){4}}
	\put(8.4,0){$c,-1$}
	\put(12.4,0){$c,1$}
	\put(10.6,-3){$b,1$}
\end{picture}
\caption{Gauss diagram and its comte}\label{fg:gdg}
\end{figure}
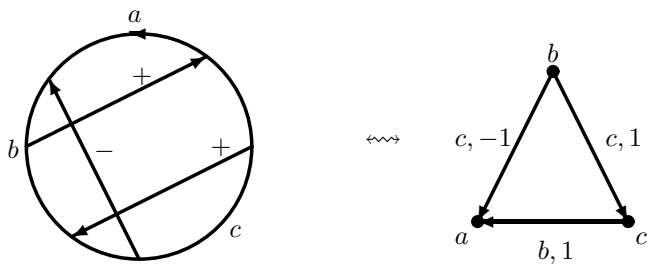

Virtual links generalize the usual (oriented) links by admitting link diagrams with, possibly,
``virtual" crossings, see \cite{kvn}.  The virtual links are the equivalence classes of
such diagrams modulo an appropriate version of the Reidemeister moves.  Another approach to
virtual links uses the classical Gauss diagrams of links. A Gauss diagram consists of several
oriented circles and arrows with distinct endpoints on the circles, see Figure~\ref{fg:gdg}
for a Gauss diagram on one circle. Each arrow should be provided with a sign $\pm 1$.
Again, there are analogs of the Reidemeister moves for Gauss diagrams and virtual links are
the equivalence classes of Gauss diagrams modulo these moves. Note that usual links can be
encoded in terms of Gauss diagrams, but not every Gauss diagram arises from a link. Thus,
virtual links can be thought as a generalization of links in which all Gauss diagrams
are allowed.

We now associate a comte with every virtual link.  Pick a Gauss diagram $G$ representing
this link.  Cut the union of the circles of $G$ at all arrowheads.  This gives a finite
number of oriented arcs $a,b,c,...$ which will be the vertices of our comte. Each arrow
$x$ of $G$ gives rise to an edge of this comte as follows.  There are two circle arcs,
say $b,c$, adjacent to the headpoint of $x$ where we choose the notation so that $b$ is
incoming and $c$ is outgoing.  We introduce an edge $\flecha bac$ if $x$ has sign $+$
and an edge $\flecha cab$ if $x$ has sign $-$, where the label $a$ is the circle arc
containing the tail of $x$.  The flow of this edge is defined to be the sign $\pm 1$ of $x$.

\begin{figure}[ht]
\setlength{\unitlength}{0.25cm}
\begin{picture}(20,8)(-4,-4)
	\allinethickness{.25mm}
	\put(0,0){\arc{8}{3.500}{4.100}}
	\put(0,0){\arc{8}{5.300}{5.900}}
	\put(0,0){\arc{8}{1.047}{2.094}}
	\put(-1,-3.8){\vector(-1,3){2.1}}
	\put(1,-3.8){\vector(1,3){2.1}}
	\put(6,0){\vector(-1,0){1.50}}
	\put(6,0){\vector(1,0){1.50}}
	\put(12,0){\arc{8}{3.500}{4.100}}
	\put(12,0){\arc{8}{5.300}{5.900}}
	\put(12,0){\arc{8}{1.047}{2.094}}
	\put(13,-3.8){\vector(-2,3){4.15}}
	\put(11,-3.8){\vector(2,3){4.15}}
\end{picture}
\caption{Move added when converting a virtual link into a comte}\label{fg:mnd}
\end{figure}
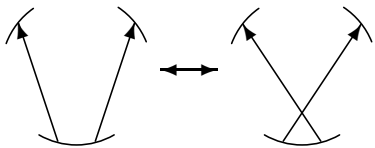

In the process of converting a virtual link into a comte, we loose some information.
Specifically, two virtual links obtained from each other by swapping two arrowtails
whenever there is no arrowhead between them, (see Figure~\ref{fg:mnd}) give rise to
the same comtes.  In the definition of virtual links in terms of link projections
with virtual crossings, there are two analogs of the third Reidemeister move for three
strands crossing at three points, two of which are virtual crossings. The two similar
moves where only one of the crossings is virtual are not allowed. If one adds these two
moves then the theory becomes empty, as any knot would become equivalent to the unknot.
By passing to the comte, we are adding one of the two forbidden moves.  Virtual links
considered modulo one of these forbidden moves are called by L.~Kauffman \emph{welded links}.

\subsection{Finite type invariants}
We outline a version of the theory of finite knot invariants for self-indexed graphs.
Note that the flows seem to play no role in this context.

For a self-indexed graph $\G=(V,E,s,t,\ell)$ and an arrow $a\in E$, we define $\G/a$ to
be the self-indexed graph $(V',E',s',t',\ell')$, where $V'=V/s(a)=t(a)$, $E'=E-\{a\}$,
and $s',t',\ell'$ are the obvious maps induced by $s,t,\ell$.  If $b$ is another arrow in
$\G$, then clearly $(\G/a)/b=(\G/b)/a$.  Thus, for any set of arrows $T\subset E$, we can
unambiguously define a self-indexed graph $\G/T$ by induction on the number of elements
of $T$.

\begin{figure}[ht]
	\setlength{\unitlength}{0.4cm}
	\begin{picture}(30,6)(-1,-1)
	\allinethickness{.5mm}
	\put(0,0){\circle*{.2}}\put(-0.6,-0.2){$e$}
	\put(0,4){\circle*{.2}}\put(-0.6,3.8){$a$}
	\put(2,2){\circle*{.2}}\put(1.2,1.8){$c$}
	\put(6,2){\circle*{.2}}\put(6.4,1.8){$d$}
	\put(8,0){\circle*{.2}}\put(8.2,-0.2){$f$}
	\put(8,4){\circle*{.2}}\put(8.2,3.8){$b$}
	\put(0,0){\vector(1,1){2}}\put(1,0.4){$d$}
	\put(2,2){\vector(-1,1){2}}\put(1,3.2){$f$}
	\drawline(2,2)(2.4,1.6)(3,1.2)(4,1)(5,1.2)(5.6,1.6)(6,2)
	\put(5.6,1.6){\vector(1,1){.4}}\put(3.8,0.4){$b$}
	\dashline{.4}(2,2)(2.4,2.4)(3,2.8)(4,3)(5,2.8)(5.6,2.4)(6,2)
	\put(5.6,2.4){\vector(1,-1){.4}}\put(3.8,3.2){$f$}
	\put(8,4){\vector(-1,-1){2}}\put(6.8,3.2){$a$}
	\put(6,2){\vector(1,-1){2}}\put(6.8,0.4){$c$}
	\put(10,2){$=$}
	\put(12,0){\circle*{.2}}\put(11.4,-0.2){$e$}
	\put(12,4){\circle*{.2}}\put(11.4,3.8){$a$}
	\put(14,2){\circle*{.2}}\put(13.2,1.8){$c$}
	\put(18,2){\circle*{.2}}\put(18.4,1.8){$d$}
	\put(20,0){\circle*{.2}}\put(20.2,-0.2){$f$}
	\put(20,4){\circle*{.2}}\put(20.2,3.8){$b$}
	\put(12,0){\vector(1,1){2}}\put(13,0.4){$d$}
	\put(14,2){\vector(-1,1){2}}\put(13,3.2){$f$}
	\drawline(14,2)(14.4,1.6)(15,1.2)(16,1)(17,1.2)(17.6,1.6)(18,2)
	\put(17.6,1.6){\vector(1,1){.4}}\put(15.8,0.4){$b$}
	\drawline(14,2)(14.4,2.4)(15,2.8)(16,3)(17,2.8)(17.6,2.4)(18,2)
	\put(17.6,2.4){\vector(1,-1){.4}}\put(15.8,3.2){$f$}
	\put(20,4){\vector(-1,-1){2}}\put(18.8,3.2){$a$}
	\put(18,2){\vector(1,-1){2}}\put(18.8,0.4){$c$}
	\put(22,2){$-$}
	\put(24,0){\circle*{.2}}\put(23.4,-0.2){$e$}
	\put(24,4){\circle*{.2}}\put(23.4,3.8){$a$}
	\put(26,2){\circle*{.2}}\put(25.2,1.8){$c$}
	\put(28,0){\circle*{.2}}\put(28.2,-0.2){$f$}
	\put(28,4){\circle*{.2}}\put(28.2,3.8){$b$}
	\put(24,0){\vector(1,1){2}}\put(25,0.4){$c$}
	\put(26,2){\vector(-1,1){2}}\put(25,3.2){$f$}
	\put(28,4){\vector(-1,-1){2}}\put(26.8,3.2){$a$}
	\put(26,2){\vector(1,-1){2}}\put(26.8,0.4){$c$}
	\spline(26,2)(27,2.6)(28,2)(27,1.4)(26,2)
	\put(26.4,1.8){\vector(-2,1){.2}}\put(28.2,1.8){$b$}
\end{picture}
\caption{A semi-virtual arrow}\label{fg:sva}
\end{figure}
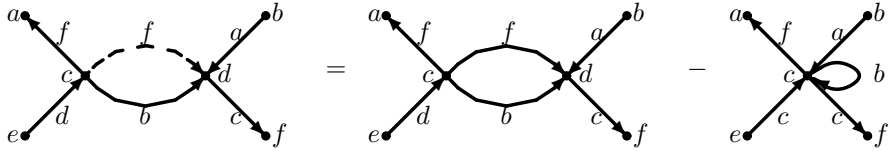

In analogy with \cite{gpv}, we define a \emph{self-indexed graph with semi-virtual arrows\/}
to be a pair $(\G,S)$, where $\G$ is a self-indexed graph and $S\subseteq E$ is a set of
arrows in $\G$. The arrows belonging to $S$ are called semi-virtual arrows of $(\G,S)$. We
draw semi-virtual arrows by dashed lines, as on the LHS of Fig. \ref{fg:sva}.
Self-indexed graphs with semi-virtual arrows yield a convenient way of encoding certain
linear combinations of self-indexed graphs. Namely, each self-indexed graph with semi-virtual
arrows $(\G,S)$ gives rise to the formal sum
$$[\G,S]=\sum_{T\subseteq S}(-1)^{\# (T)}\;\G/T \in \ZZ\GG$$
where $\GG$ is the set of isomorphism classes of self-indexed graphs and $\ZZ\GG$ is the
abelian group freely generated by $\GG$.

Let $\nu$ be a map from $ \GG$ to an abelian group $A$. We say that $\nu$ is \emph{an
invariant of finite type\/} if there exists an integer $n\ge 0$ such that the linear extension
$\ZZ \GG\to A$ of $\nu$ vanishes on $[\G,S]$ for all self-indexed graphs $(\G,S)$ with $n$
semi-virtual arrows.  The minimal $n$ with that property is called the \emph{degree\/}
of $\nu$. For instance, an invariant of degree $1$ assigns the same element of $A$ to
all self-indexed graphs.  Of course, in our context we are interested only in those $\nu$
which are preserved under the transformations R0-R3.

\subsection{Remark}
Other objects from knot theory have their counterparts in the world of self-indexed graphs
and comtes.  This includes knotted graphs in $\RR^3$, braids, and tangles.  It would be
interesting to reformulate further classical knot invariants, in particular the Conway
and Jones polynomials, in terms of comtes.  The authors plan to consider some of these
questions in another publication.

\section{Homology of self-indexed q-graphs}
In this section we generalize the quandle cocycle invariants of comtes, replacing quandles
by so-called self-indexed q-graphs.  We also introduce a homology theory for the self-indexed
q-graphs generalizing the homology of quandles.

\subsection{Homology of self-indexed graphs}
For each $m\ge 1$, consider the graph $y_m$ obtained as the 1-skeleton of the unit cube
in $\RR^{m-1}$. More precisely, let $\{e_1,\ldots,e_{m-1}\}$ be the canonical basis of
$\RR^{m-1}$. The vertices of $y_m$ are the points $(i_1,\ldots,i_{m-1})\in\RR^{m-1}$ such
that $i_j\in\{0,1\}$.  There is an arrow in $y_m$ pointing form a vertex $v$ to a vertex
$w$ iff $w-v=e_s$ for some $s$.
We now label the vertices and arrows of $y_m$ by sequences of positive integers.  The labels
are defined by induction on $m$. The only vertex of $y_1$ has label $1$. Assume that $y_m$
is labeled. Consider the intersection of $y_{m+1}$ with the hyperplane $\{x_1=0\}$ and
copy to it the labels from $y_m$ with a shift: if a vertex or an arrow of $y_m$ has a label
$j_1,j_2,\ldots, j_r$, then its copy in the hyperplane has the label $ j_1+1,j_2+1,\ldots,
j_r+1$. All the arrows of $y_{m+1}$ parallel to $e_1$ are labeled with $1$. The labels
of vertices and arrows of $y_{m+1}$ lying in the hyperplane $\{x_1=1\}$ are those of the
hyperplane $\{x_1=0\}$ with a $1$ added at the beginning: if a vertex or an arrow $v$
of $\{x_1=0\}$ has a label $j_1,j_2,\ldots, j_r$, then its parallel translation $v+e_1$
has the label $1,j_1,j_2,\ldots, j_r$.

We provide the disjoint union $Y_n=\amalg_{m=1}^n y_m$ with the structure of a self-indexed
graph: the map from the set of arrows of $Y_n$ to the set of its vertices is defined by
assigning to each arrow the only vertex with the same label.
We draw $y_m$ for $m=1,2,3,4$ in Figure~\ref{fg:yn}. The whole figure is then a drawing of $Y_4$.
\begin{figure}[ht]
	\setlength{\unitlength}{0.56cm}
	\begin{picture}(26,8)(-19,-1)
	\allinethickness{.5mm}
	\put(-18.2,-1){$y_1$}
	\put(-18,3){\circle*{.25}}
	\put(-18.2,2.2){$\mathbf{1}$}
	\put(-14.4,-1){$y_2$}
	\put(-16,3){\circle*{.25}}
	\put(-12,3){\circle*{.25}}
	\put(-16,3){\vector(1,0){4}}
	\put(-16.2,2.2){$\mathbf{2}$}
	\put(-12.3,2.2){$\mathbf{12}$}
	\put(-14.2,3.2){$\scriptstyle 1$}
	\put(-6,-1){$y_3$}
	\put(-8,1){\circle*{.25}}
	\put(-4,1){\circle*{.25}}
	\put(-8,5){\circle*{.25}}
	\put(-4,5){\circle*{.25}}
	\put(-8,1){\vector(1,0){4}}
	\put(-8,1){\vector(0,1){4}}
	\put(-8,5){\vector(1,0){4}}
	\put(-4,1){\vector(0,1){4}}
	\put(-8.5,0.5){$\mathbf{3}$}
	\put(-8.8,5){$\mathbf{23}$}
	\put(-3.7,0.5){$\mathbf{13}$}
	\put(-3.8,5){$\mathbf{123}$}
	\put(-7.8,3){$\scriptstyle 2$}
	\put(-3.8,3){$\scriptstyle 12$}
	\put(-6.2,1.2){$\scriptstyle 1$}
	\put(-6.2,5.2){$\scriptstyle 1$}
	\put(2,-1){$y_4$}
	\put(0,0){\circle*{.25}}
	\put(4,0){\circle*{.25}}
	\put(0,4){\circle*{.25}}
	\put(4,4){\circle*{.25}}
	\put(2,2){\circle*{.25}}
	\put(6,2){\circle*{.25}}
	\put(2,6){\circle*{.25}}
	\put(6,6){\circle*{.25}}
	\put(2,2){\line(1,0){1.8}}\put(4.2,2){\vector(1,0){1.8}}
	\put(2,2){\line(0,1){1.8}}\put(2,4.2){\vector(0,1){1.8}}
	\put(0,0){\vector(1,1){2}}
	\put(6,2){\vector(0,1){4}}
	\put(4,0){\vector(1,1){2}}
	\put(0,0){\vector(1,0){4}}
	\put(0,0){\vector(0,1){4}}
	\put(2,6){\vector(1,0){4}}
	\put(0,4){\vector(1,1){2}}
	\put(4,4){\vector(1,1){2}}
	\put(0,4){\vector(1,0){4}}
	\put(4,0){\vector(0,1){4}}
	\put(-0.5,-0.5){$\mathbf{4}$}
	\put(1.0,1.8){$\mathbf{24}$}
	\put(-0.9,3.8){$\mathbf{34}$}
	\put(0.6,5.8){$\mathbf{234}$}
	\put(4.1,-0.5){$\mathbf{14}$}
	\put(6.2,1.8){$\mathbf{124}$}
	\put(4.2,3.8){$\mathbf{134}$}
	\put(6.3,5.8){$\mathbf{1234}$}
	\put(1.3,0.9){$\scriptstyle 2$}
	\put(5.3,0.9){$\scriptstyle 12$}
	\put(1.3,4.9){$\scriptstyle 2$}
	\put(5.3,4.9){$\scriptstyle 12$}
	\put(3,2.2){$\scriptstyle 1$}
	\put(3,6.2){$\scriptstyle 1$}
	\put(1,0.2){$\scriptstyle 1$}
	\put(1,4.2){$\scriptstyle 1$}
	\put(2.1,3){$\scriptstyle 23$}
	\put(6.1,3){$\scriptstyle 123$}
	\put(0.2,1.2){$\scriptstyle 3$}
	\put(4.2,1.2){$\scriptstyle 13$}
	\end{picture}
\caption{Graph $Y_4=y_1\cup y_2\cup y_3\cup y_4$.}\label{fg:yn}
\end{figure}

We consider $2n$ embeddings of $Y_n$ in $Y_{n+1}$. They are the only homomorphisms
compatible with the embeddings of $y_n$ to $y_{n+1}$ given by
\begin{align*}
D_s^0(i_1,\ldots,i_{n-1}) &= (i_1,\ldots,i_{s-1},0,i_s,\ldots,i_{n-1}) \\
D_s^1(i_1,\ldots,i_{n-1}) &= (i_1,\ldots,i_{s-1},1,i_s,\ldots,i_{n-1}),
\end{align*}
where $s=1,2,\ldots,n$. These homomorphisms are the \emph{faces} of $Y_{n+1}$.
We can give an explicit description by considering the labels in the free
quandle of the set $\{1,2,\ldots,n+1\}$. The label $j_1,\ldots,j_r$ corresponds
in this description to $j_1\trid(j_2\trid(\cdots\trid j_r))$. Then
\begin{align*}
D_s^0(j_1\trid(\cdots\trid (j_{r-1}\trid j_r))) &= i_1\trid(\cdots\trid (i_{r-1}\trid i_r)),\qquad
	\text{with }i_t=\begin{cases}j_t+1&\text{if }j_t\ge s,\\j_t&\text{if }j_t<s,\end{cases} \\
D_s^1(j_1\trid(\cdots\trid (j_{r-1}\trid j_r))) &= i_1\trid(\cdots\trid (i_{r-1}\trid i_r)),\qquad
	\text{with }i_t=\begin{cases}s\trid (j_t+1)&\text{if }j_t\ge s,\\j_t&\text{if }j_t<s.\end{cases} \\
\end{align*}
For instance, let us consider the embedding $D_2^1:Y_3\to Y_4$. We get
$D_2^1(3)=2\trid 4$, $D_2^1(1\trid 3)=1\trid (2\trid 4)$,
$D_2^1(2\trid 3)=(2\trid 3)\trid (2\trid 4)=2\trid (3\trid 4)$,
$D_2^1(1\trid(2\trid 3))=1\trid ((2\trid 3)\trid(2\trid 4))=1\trid(2\trid(3\trid 4))$.
Thus, $D_2^1$ sends $y_3$ to the ``rear" face of $y_4$. Analogously, it sends $y_2$ to the
``top" face of $y_3$, with vertices $23,123$; and it sends $y_1$ to itself.

\bigskip
We now define a homology theory for self-indexed graphs.
\begin{defin}
Let $\G=(V,E,s,t,\ell)$ be a self-indexed graph. Let $C_n(\G)$ be the free abelian
group generated by all homomorphisms $f: Y_n\to\G$.  Define the boundary map
$\partial :C_{n}(\G)\to C_{n-1}(\G)$ by
$$\partial (f)=\sum_{s=1}^{n-1}(-1)^s(f D_s^0-f D_s^1) .$$
This gives a chain complex $C_\bullet(\G) $.  For an abelian group $A$, set
$$C_\bullet(\G,A)=C_\bullet(\G)\otimes A,\quad
	C^\bullet(\G,A)=\Hom(C_\bullet(\G),A).$$
Consequently, we obtain homology and cohomology theories from these chain complexes.
\end{defin}

We warn that our homology has nothing to do with Kontsevich's graph homology theory.
In his theory one fixes a species (or an operad) and graphs are used
to define a basis for a chain complex. Here we fix a self-indexed graph and the basis is
given by certain homomorphisms.

\subsection{Self-indexed r-graphs and q-graphs}
A self-indexed graph $\G=(V,E,s,t,\ell)$ is an \emph{r-graph\/} if
\begin{enumerate}
\item two different arrows pointing out from the same vertex have different labels, and
\item two different arrows pointing into the same vertex have different labels.
\end{enumerate}
We say that $\G$ is a \emph{q-graph\/} if furthermore
\begin{enumerate}\addtocounter{enumi}{2}
\item each vertex $a\in V$ has an arrow $\flecha aaa$ with label, source and target $a$.
\end{enumerate}
If $\flecha bac$ is an arrow in a self-indexed r-graph, then we write $c=a\cdot b$.

For instance, the self-indexed graph of a rack is a self-indexed r-graph. The self-indexed
graph of a quandle is a self-indexed q-graph.  However, there are many more self-indexed
r-graphs and self-indexed q-graphs than racks and quandles (see for instance Example
\ref{ex:hoc}).

If $\G$ is a self-indexed r-graph, then a homomorphism $f:Y_n\to\G$ is uniquely determined
by $f(1),\ldots,f(n)$ (though it is not true that any sequence $f(1),\ldots,f(n)$ gives rise
to a homomorphism). In this case, we denote $f$ by $<f(1),\ldots,f(n)>$.  The definition
of the boundary $C_n(\G)\to C_{n-1}(\G)$ can be re-written as
$$\partial (<a_1,\ldots,a_{n}>)
	= \sum_{s=1}^{n-1} (-1)^s (<a_1,\ldots,a_{s-1},a_{s+1},\ldots,a_{n}>
		-<a_1,\ldots,a_{s-1},a_s\cdot a_{s+1},\ldots,a_s\cdot a_{n}> ).$$

If $\G$ is a self-indexed q-graph, we define $C^Q_n(\G)$ to be $C_n(\G)$ quotiented out
by the morphisms $<a_1,\ldots,a_n>$ such that $a_i=a_{i+1}$ for some $i$.  This gives a
quotient chain complex $C_\bullet^Q(\G)$.  Set
$$
C_\bullet^Q(\G,A)=C_\bullet^Q(\G)\otimes A,\quad
C^\bullet_Q(\G,A)=\Hom(C_\bullet^Q(\G),A).$$ We obtain (co)homology theories from these
chain complexes.  We refer to cycles and cocycles in $C_\bullet^Q$, $C^\bullet_Q$ as
q-cycles and q-cocycles.

It is easy to see that if $\G$ is the self-indexed graph of a rack then any $n$-tuple
$<a_1,\ldots,a_n>$ gives rise to a homomorphism and the chain complex $C_\bullet(\G)$
coincides with the usual chain complex of the rack.  Also, if $\G$ is the self-indexed graph
of a quandle, then $C_\bullet^Q(\G)$ coincides with the usual chain complex of the quandle
(cf. \cite{cjkls}).

Our (co)homology of $\G$ in general is not invariant under the moves R0--R3 on $\G$. For
instance, let $\G_2$ and $\G_3$ be the self-indexed graphs respectively in the middle and
RHS of Figure~\ref{fg:gah} below. If $\bar\G_2$ and $\bar\G_3$ denote the self-indexed graphs
obtained by adding the arrows $\flecha iii$ for $i=a,b,c$, then we have $H_3(\bar\G_2)=\ZZ^4$,
while $H_3(\bar\G_3)=\ZZ^5$. However, these two self-indexed graphs are isotopic.

\medskip
We classified self-indexed r-graphs with $3$ vertices. Modulo isomorphism, there are $6663$
such self-indexed graphs, among which $70$ are self-indexed q-graphs. We computed their
homology with integer coefficients up to degree $5$.  There are $280$ different values
for self-indexed r-graphs and $28$ for self-indexed q-graphs. We just give an example here.
\begin{examp}\label{ex:hoc}
Let $V=\{a,b,c\}$ and consider arrows $\flecha iii$ for $i=a,b,c$, and $\flecha bab$,
$\flecha cac$, $\flecha abc$, $\flecha cba$, $\flecha acb$.  The homology up to degree $5$
of this self-indexed graph is $H_1 = \ZZ$, $H_2 = \ZZ^2$, $H_3 = \ZZ^4$, $H_4 = \ZZ^7$,
$H_5 = \ZZ^{11}$, from where it seems that the $n$-th Betti number of this self-indexed
graph is $ \frac 12n(n-1)+1$.  Note that there is no rack with such homology, as Betti
numbers of racks grow exponentially (see \cite{eg}).
\end{examp}

\subsection{Quandle cocycle invariants re-examined}\label{ss:vgd}
Fix a self-indexed graph $\G$. Any $n$-chain $I\in C_n(\G)$ uniquely expands as
$$I=\sum_{\tau:Y_n\to\G} I_{\tau}\tau.$$
For an $n$-cochain $f\in C^n(\G,A)$ with values in an abelian group $A$, set
$$\int_{\G}I\,f=\sum_{\tau:Y_n\to\G} I_\tau f(\tau).$$
Similarly, for a homomorphism of self-indexed graphs $\sigma:\G'\to\G$, an $n$-chain $I$
on $\G'$ and an $A$-valued $n$-cochain $f$ on $\G$, set
$$
\int_\sigma I\,f=\int_{\G'}I\,(\sigma^*f)= \sum_{\tau:Y_n\to\G'}I_\tau f(\sigma\circ \tau).
$$
It is clear that if $J$ is an $(n+1)$-chain on $\G'$, then
\begin{equation}\label{eq:idd}
\int_\sigma (\partial J)\,f = \int_\sigma J\,(\partial^*f).
\end{equation}
Therefore $\int_\sigma$ induces a bilinear pairing $H_n(\G')\otimes H^n(\G;A)\to A$.

We can use this formalism to define a \emph{state sum\/} on a pair (a self-indexed graph
$\G'$, a chain $I\in C_n(\G')$). For a cochain $f\in C^n(\G,A)$, set
\begin{equation}\label{df:qigg}
\invar If=\sum_{\substack{\text{homomorphisms} \\ \sigma:\G'\to\G}}
	\int_\sigma I\,f\ \in\ \ZZ A.
\end{equation}
The sum here is the sum in the ring $\ZZ A$, while the integral on the RHS is given by
sums in $A$. If $I$ is a cycle and $f$ is a cocycle, then $\invar If$ depends only on their
(co)homology classes, since by \eqref{eq:idd},
$$\invar {(\partial J)}f=\invar J{(\partial^*f)}.$$
This relation holds also if $\G$ is a q-graph and $f\in C^\bullet_Q(\G;A)$.

\begin{corol}\label{co:pdoh}
The state sum $\invar If$ defines a pairing $\invar{}{}:H_n(\G')\times H^n(\G;A)\to\ZZ A$
(which, we warn, is not bilinear).  If $\G$ is a q-graph, then this state sum defines also
a pairing $\invar{}{}:H_n(\G')\times H^n_Q(\G;A)\to\ZZ A$.
\hfill\qed
\end{corol}

\subsection{Degree $2$}
We focus now on the case $n=2$. Fix an abelian group $A$.  Note that a flow on a self-indexed
graph $\G'$ is the same thing as a cycle in $C_2(\G')$.  If a self-indexed q-graph $\G$
is derived from a quandle, then, as remarked above, a cocycle in $C^2_Q(\G;A)$ is just a
$2$-cocycle of this quandle with values in $A$.  Thus, let $(\G',I)$ be a comte, let $X$
be a finite quandle, let $\G$ be its self-indexed q-graph and let $f\in Z^2_Q(\G;A)$. It
follows from definitions that under these assumptions
$$\invar If=\Phi(\G',I,\G,f).$$

\begin{propo}
Let $\G$ be a self-indexed q-graph and $f\in Z^2_Q(\G;A)$ a cocycle.  The state sum
$\invar If$ does not change under the moves R1, R2 on the comte $(\G',I)$. If $\G$ is
derived from a quandle, then $\invar If$ does not change under all the moves R0--R3 on
$(\G',I)$. Also $\invar If=\invar I{f'}$ if $f$ is cohomologous to $f'$.
\end{propo}
\begin{proof}
The invariance of $\invar If$ under the moves R2(a) and R2(b) follow from $\G$ being a
self-indexed r-graph: if the LHS in Fig.~\ref{fg:r2a} is part of $\G'$, then any homomorphism
$ \G'\to\G$ must send $s$ and $t$ to the same vertex. The invariance of $\invar If$ under
R1 follows from $\G$ being a self-indexed q-graph: if an arrow $e$ of $\G'$ is labeled by
one of its endpoints, then any homomorphism $ \G'\to\G$ must send both endpoints to the
same vertex. Also, since $f$ is a q-cocycle, any arrow labeled with its source gives no
contribution to the state sum invariant.

Now, suppose $\G$ is the self-indexed graph of a quandle.  The invariance of $\invar If$
under R0 follows since the flow on the deleted arrow is $0$ (note that any homomorphism
defined on the self-indexed graph on the RHS of Fig.~\ref{fg:r0} can be extended in a
unique way to the self-indexed graph on the LHS).  For the move R3(a), let $\G'$ contain
a subgraph as on the RHS of Fig. \ref{fg:r3a}. The condition $x\trid (y\trid z)=(x\trid
y)\trid (x\trid z)$ implies that any graph homomorphism $\sigma:\G'\to\G$ uniquely extends
to a homomorphism from the self-indexed graph on the LHS to $\G$ (and analogously for the
other edges in the square).  Since the flow on the edge is assumed to be $0$, the state sum
does not change.  As for R3(b), this is an immediate consequence of the $2$-cocycle condition.

The last assertion is a consequence of Corollary \ref{co:pdoh}.
\end{proof}

\begin{remar}
The requirement on $\G$ to be a quandle is necessary to have a bijection
of the sets of graph homomorphisms for the self-indexed graphs on both sides
of Fig.~\ref{fg:r0} and Fig.~\ref{fg:r3a}. In order for $\invar If$ to be invariant
under R3, it would be enough to assume that $\G$ has the property that
whenever one has a subgraph as on the RHS of Fig. \ref{fg:r3a}, the fourth arrow
exists in $\G$, as on the LHS of Fig.~\ref{fg:r3a}. However, in this case $\invar If$
would not be necessarily invariant under R0. In particular, considered for knots, it
would not be invariant under the second and the third Reidemeister moves.
\end{remar}

\subsection{Examples}
The smallest indecomposable quandle with $H^2_Q\neq 0$ can be identified with the vertices
of a tetrahedron.  We denote its elements $\{0,1,2,3\}$. Each vertex $i$ acts by a rotation
of the tetrahedron by an angle of $\frac{2\pi}3$ fixing $i$, (see Figure~\ref{fg:tetra},
where the action of $0$ is drawn).

\begin{figure}[ht]
\setlength{\unitlength}{1cm}
\begin{picture}(5.5,2.5)(-.5,.5)
	\allinethickness{.5mm}
	\put(1,1){\circle*{.1}}
	\put(3,1){\circle*{.1}}
	\put(3,2){\circle*{.1}}
	\put(2,3){\circle*{.1}}
	\put(0.7,0.8){$1$}
	\put(3.2,0.8){$2$}
	\put(3.2,1.9){$3$}
	\put(1.9,3.1){$0$}
	\path(1,1)(3,1)(2,3)(1,1)
	\path(3,1)(3,2)(2,3)
	\dottedline[.]{.15}(1,1)(3,2)
	\spline(1.2,0.8)(2,0.6)(2.8,0.8)
	\spline(3.2,1.2)(3.4,1.5)(3.2,1.8)
	\put(2.8,0.8){\vector(2,1){.1}}
	\put(3.2,1.8){\vector(-1,2){.1}}
	\put(1.8,0.3){$0\trid ?$}
	\put(3.7,1.3){$0\trid ?$}
\end{picture}
\caption{}\label{fg:tetra}
\end{figure}

Alternatively, one can think of this quandle as the affine (=Alexander) quandle over the
field with $4$ elements $\FF_4$, and the automorphism determined by $\omega\in\FF_4-\{0,1\}$.
Specifically, $x\trid y=(1-\omega)x+\omega y$. An isomorphism $f:\text{tetrahedron}\to\FF_4$
is given by $f(0)=0$, $f(1)=1$ $f(2)=\omega$, $f(3)=1+\omega$ (see \cite{ag}). This quandle
has a nontrivial $2$-cocycle with values in the cyclic group with two elements, $C_2$. Let
$\sigma$ be the generator of $C_2$; we can write the cocycle as
$$f(<x,y>) =\sigma^{1-\delta(xy(x-y))},$$
where $\delta(z)=1$ if $z=0$ and $\delta(z)=0$ otherwise.

We compute the invariant of comtes derived from this 2-cocycle for the comtes $\G_1,\G_2,\G_3$
in Figure~\ref{fg:gah}.  It is easy to see that $\invar {I_1}f=4+12\sigma$ and
$\invar {I_2}f=\invar {I_3}f=4$.  Indeed, $\G_2$ and $\G_3$ are related by a sequence of moves
R1, R3(a), R3(b), R3(a), R1.

\begin{figure}[ht]
\setlength{\unitlength}{.5cm}
\begin{picture}(26,7)(-1,1)
	\allinethickness{.5mm}
	\put(1,4){\circle*{.25}}
	\put(3,8){\circle*{.25}}
	\put(5,4){\circle*{.25}}
	\put(0.5,3.5){$a$}
	\put(2.85,8.25){$b$}
	\put(5.1,3.5){$c$}
	\put(1,4){\vector(1,2){2}}
	\put(3,8){\vector(1,-2){2}}
	\put(5,4){\vector(-1,0){4}}
	\put(1,5.75){$c,1$}
	\put(4.25,5.75){$a,1$}
	\put(2.5,4.25){$b,1$}
	\put(2.25,1){$\G_1,I_1$}
	\put(9.5,4){\circle*{.25}}
	\put(11.5,8){\circle*{.25}}
	\put(13.5,4){\circle*{.25}}
	\put(9,3.5){$a$}
	\put(11.35,8.25){$b$}
	\put(13.6,3.5){$c$}
	\put(9.5,4){\vector(1,2){2}}
	\put(11.5,8){\vector(1,-2){2}}
	\put(13.5,4){\vector(-1,0){4}}
	\put(9.5,5.75){$c,1$}
	\put(12.75,5.75){$a,1$}
	\put(11,4.25){$b,1$}
	\spline(9.5,4)(11.5,3)(13.5,4)
	\put(13.4,3.9){\vector(2,1){.1}}
	\put(11,2.35){$b,0$}
	\put(10.75,1){$\G_2,I_2$}
	\put(18,4){\circle*{.25}}
	\put(20,8){\circle*{.25}}
	\put(22,4){\circle*{.25}}
	\put(17.5,3.5){$a$}
	\put(19.85,8.25){$b$}
	\put(22.1,3.5){$c$}
	\put(18,4){\vector(1,2){2}}
	\put(22,4){\vector(-1,0){4}}
	\put(19.25,5.75){$c,1$}
	\put(17,5.75){$c,1$}
	\put(19.5,4.25){$b,0$}
	\spline(18,4)(20,3)(22,4)
	\put(21.9,3.9){\vector(3,1){.1}}
	\spline(20,8)(18,6)(18,4.1)
	\put(18,4.2){\vector(0,-1){.1}}
	\put(19.5,2.35){$b,0$}
	\put(19.25,1){$\G_3,I_3$}
\end{picture}
\caption{}\label{fg:gah}
\end{figure}

\subsection{Acknowledgements}
This work was written at IRMA, Strasbourg, while M.G. visited V.T.
M.G. is deeply indebted to the whole IRMA for its warm hospitality;
his thanks go also to L. Kauffman for correspondence about welded knots.
V.T. would like to thank
S. Matveev for useful correspondence on the quandles of links.

\end{document}